\documentclass[oneside,english,12pt,reqno]{amsart}
\usepackage[colorlinks=true,linkcolor=blue,citecolor=magenta]{hyperref}
\usepackage{comment}

\allowdisplaybreaks

\usepackage[mathscr]{eucal}
\usepackage{typearea}
\usepackage{charter}
\usepackage{graphicx}

\usepackage{units}
\usepackage{mathrsfs}
\usepackage{amstext}
\usepackage{amsthm}
\usepackage{esint}
\usepackage{amsmath}
\usepackage{amssymb}
\usepackage{soul}
\usepackage{cancel}
\newcommand{\R}{\mathbb{R}}
\newcommand{\pa}{\partial}
\newcommand{\eps}{\varepsilon}
\usepackage{a4wide}
\usepackage{enumitem}
\newcommand{\bra}[1]{\left(#1\right)}

\newcommand{\na}{\nabla}

\renewcommand{\L}{\mathscr{L}}
\renewcommand{\H}{\mathscr{H}}

\newcommand{\ue}{u^{\eps}}

\newcommand{\sumi}{\sum_{i=1}^{m}}

\makeatletter
\numberwithin{equation}{section}
\numberwithin{figure}{section}
\newtheorem{corollary}{Corollary}[section]
\newtheorem{theorem}{Theorem}[section]

\newtheorem{remark}{Remark}[section]
\newtheorem{definition}{Definition}[section]
  
\makeatother

\usepackage{babel}

\title[Global existence for reaction-diffusion systems on multiple domains]{Global existence for reaction-diffusion systems on multiple domains}
\author[W.E. Fitzgibbon]{William E. Fitzgibbon}
\address{William E. Fitzgibbon \hfill\break
	Department of Mathematics, 
	University of Houston, Houston, Texas 77004, USA}
\email{fitz@uh.edu}
\author[J.J. Morgan]{Jeff Morgan}
\address{Jeff Morgan \hfill\break
	Department of Mathematics, 
	University of Houston, Houston, Texas 77004, USA}
\email{jmorgan@uh.edu}
\author[J.M. Ryan]{John Maurice-Car Ryan}
\address{John Ryan \hfill\break
	Department of Mathematics, 
	University of Houston, Houston, Texas 77004, USA}
\email{jmryan2@central.uh.edu}

\begin{document}
\begin{abstract}
	We study the global existence of solutions reaction-diffusion systems with control of mass on multiple domains. Some of these domains overlap, and as a result, an unknown defined on one subdomain can impact another unknown defined on a different domain that intersects with the first. Our results extend those in \cite{fmty}.
\end{abstract}
\maketitle
\tableofcontents

\section{Problem Setting}

\subsection{Introduction}
This work is concerned with reaction-diffusion systems that are defined on a sequence of spatially bounded non-coincident spatial subdomains $\Omega_1,...,\Omega_N\subset\R^n$. The systems allow for discontinuity in the coefficients of the differential operators and in the components of the reaction vector fields, as well as interaction of species on overlapping subdomains. The vector field is required to satisfy a quasi-positivity condition to preserve nonnegativity, and also satisfy properties that help preserve total mass/concentration. Our concern is the establishment of a priori bounds and global existence of sup norm bounded weak solutions of these systems. There has been a wealth of information for systems of this type with smooth coefficients on the differential operators and continuous reaction vector fields in the case that N=1.  In this regard, we cite \cite{desvillettes2007global, souplet2018global,pierre2017dissipative, canizo2014improved,caputo2009global, morgan1990boundedness,hollis1987global,fitzgibbon1997stability, pierre2000blowup, pierre2010global, fellner2021uniform, fellner2020global, desvillettes2007global, morgan2020boundedness, morgan1989global}. Global existence and boundedness of weak solutions in the case of differential operators with  discontinuous coefficients and discontinuous reaction vector fields have recently appeared in \cite{fmty}. There have been a few results that extend the results on single domains to results coupled across multiple domains, and we list \cite{FLMM} and \cite{FLM1}. The work at hand differs from \cite{FLMM} and \cite{FLM1} by virtue of the fact that the diffusion and the reaction vector fields can be more complex. The present work is an extension of work in \cite{fmty} for single domains with $L^\infty$ diffusion.

Problems of this type can arise in the modelling of biological
systems, and have been studied as mathematical
models. For example, one such system which is analyzed in \cite{FLM1} models
the interaction of two hosts and a vector population, where a disease is transmitted in a criss-cross fashion
from one host through a vector to another host. It is assumed that
the disease is benign for one host and lethal to the other. 

In order to provide a more complete example of what we have in mind we provide three examples. Two of these examples are given below, and revisited in Section 4, and the third example is introduced and discussed in Section 4. The first example concerns the cross species spatial transmission of an infectious disease, and the second example concerns a hypothetical interaction of three species living on two overlapping domains. In the first case, we consider an infectious disease that can be transmitted across multiple species and multiple habitats. These are a major concern for animal husbandry, wildlife management, and human health \cite{multiauthor}. A species occupying a given habit may contract the disease from a second species occupying an overlapping habit and via dispersion transmit the disease to a third species whose habitat also overlaps the habitat of the first species. In the second setting, species $A$, $B$ and $C$ interact through a reaction of the form $A+B\leftrightarrow C$ on overlapping domains $\Omega_1$ and $\Omega_2$. Species $A$ lives on $\Omega_1$, while species $B$ and $C$ live on $\Omega_2$. 
 
\subsection{Two Illustrative Examples}

Consider a spatially distributed population. The dispersion of
the population is modeled by Fickian diffusion. In this model
there are three populations confined to separate habitats $\Omega_1$, $\Omega_2$ and $\Omega_3$, such that $\Omega_1 \cap \Omega_2 \ne \emptyset$, $\Omega_2\cap\Omega_3\ne\emptyset$ and $\Omega_1\cap\Omega_3=\emptyset$. The possibility of physically separated habitats for the
vulnerable and resistant hosts are allowed, each of which
intersects with the domain of the vector.

Suppose $k_{1},$ $k_{2},$ $k_{3}$ and $k_{4}$ are nonnegative functions, 
$\lambda _{1},$ $\lambda _{2}$ and $\lambda _{3}$ are positive
constants. Furthermore, the supports of $k_{1}$ and $k_{2}$ are
contained in the intersection of $\Omega _{1}$ and $\Omega _{2}$, respectively, and
the supports of \ $k_{3}$ and $k_{4}$ are contained in the
intersection of $\Omega _{2}$ and $\Omega _{3}$, respectively. Finally, for each
$i=1,2,..6$, $d_{i}$ is a positive bounded function that is bounded away from $0$, and for each $j=1,2,3$, $\lambda _{j}$ is a positive constant. 

\begin{equation}\label{eqex1}
\left\{
\begin{array}{cc}
\left(
\begin{array}{c}
\phi _{t}=\nabla\cdot\left(d_{1}\nabla \phi\right) -k_{1}(x)\phi \beta +\lambda _{1}\psi  \\
\psi _{t}=\nabla\cdot\left(d_{2}\nabla \psi\right) +k_{1}(x)\phi \beta -\lambda _{1}\psi
\end{array}
\mbox{ for }x\in \Omega _{1},t>0\right) & \text{host 1} \\
&\\
\left(
\begin{array}{c}
\alpha _{t}=\nabla\cdot\left(d_{3}\nabla \alpha\right) -k_{2}(x)\alpha \psi
-k_{3}(x)\alpha
v+\lambda _{2}\beta  \\
\beta _{t}=\nabla\cdot\left(d_{4}\nabla \beta\right) +k_{2}(x)\alpha \psi +k_{3}(x)\alpha
v-\lambda _{2}\beta
\end{array}
\mbox{ for }x\in \Omega _{2},t>0\right) & \text{vector} \\
&\\
\left(
\begin{array}{c}
v_{t}=\nabla\cdot\left(d_{5}\nabla v\right)-k_{4}(x)v\beta  \\
w_{t}=\nabla\cdot\left(d_{6}\nabla w\right)+k_{4}(x)v\beta -\lambda _{3}w
\end{array}
\mbox{ for }x\in \Omega _{3},t>0\right) &  \text{host 2}
\end{array}
\right.
\end{equation}
\medskip 

\begin{figure}[h!]
  \includegraphics[width=\linewidth]{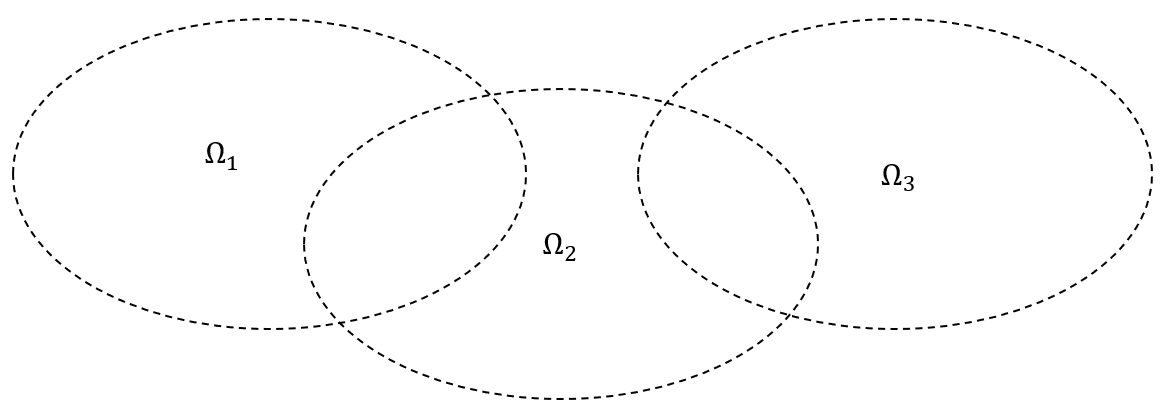}
  \caption{}
  \label{fig:domains3}
\end{figure}

We impose homogeneous Neumann boundary conditions on each domain $\Omega _{1},\Omega _{2},$ and $\Omega _{3}$.

\begin{equation}\label{eqex1bc}
\left\{
\begin{array}{cc}
\partial \phi /\partial \eta =\partial \psi /\partial \eta =0 & \text{for }
x\in \partial \Omega _{1},t>0, \\
\partial \alpha /\partial \eta =\partial \beta /\partial \eta =0 & \text{for
}x\in \partial \Omega _{2},t>0, \\
\partial v/\partial \eta =\partial w/\partial \eta =0 & \text{for }x\in
\partial \Omega _{3},t>0,
\end{array}
\right.
\end{equation}

\medskip\noindent Finally, we specify continuous nonnegative initial data.

\begin{equation}\label{eqex1ic}
\left\{
\begin{array}{cc}
\phi (0,x)=\phi _{0}(x),\quad \psi (0,x)=\psi _{0}(x) & \mbox{for }x\in
{\Omega} _{1}\\
\alpha (0,x)=\alpha _{0}(x),\quad \beta (0,x)=\beta _{0}(x) & \mbox{for }
x\in {\Omega} _{2} \\
v(0,x)=v_{0}(x),\quad w(0,x)=w_{0}(x) & \mbox{for }x\in
{\Omega} _{3}
\end{array}
\right.
\end{equation}

\medskip It can be shown that the system above preserves the nonnegativity of the initial data. In addition, on $\Omega _{1}$ the
vector field
\begin{equation}
\left(
\begin{array}{c}
-k_{1}(x)\phi \beta +\lambda _{1}\psi  \\
+k_{1}(x)\phi \beta -\lambda _{1}\psi
\end{array}
\right)   \label{eq116}
\end{equation}
has a first component that is bounded above by a linear expression, and the components that clearly sum to zero. Similarly, on $\Omega
_{2}$ the vector field
\begin{equation}
\left(
\begin{array}{c}
-k_{2}(x)\alpha \psi -k_{3}(x)\alpha v+\lambda _{2}\beta  \\
k_{2}(x)\alpha \psi +k_{3}(x)\alpha v-\lambda _{2}\beta
\end{array}
\right)   \label{eq117}
\end{equation}
has a first component that is bounded above by a linear expression, and also sums to zero. The same mechanism can
be seen on $\Omega _{3}$ since the function
\begin{equation}
\left(
\begin{array}{c}
-k_{4}(x)v\beta  \\
k_{4}(x)v\beta -\lambda _{3}w
\end{array}
\right)   \label{eq1170}
\end{equation}
has a first component that is bounded above by a linear expression, and a sum that is less than or equal to zero. We will apply our results to system (\ref{eqex1})-(\ref{eqex1ic}), and a slightly more complex extension, in Section 5. 

\medskip The second example is easier to state. Here, species $A$, $B$ and $C$ interact through a reaction of the form $A+B\leftrightarrow C$ on overlapping domains $\Omega_1$ and $\Omega_2$. Species $A$ occupies $\Omega_1$, while species $B$ and $C$ occupy $\Omega_2$. If we define $k(x)=\chi_{\Omega_1\cap\Omega_2}(x)$ (the characteristic function on $\Omega_1\cap\Omega_2$), and use $u_1(t,x)$, $u_2(t,x)$ and $u_3(t,x)$ to denote the concentration densities of $A$, $B$ and $C$, then a possible model is given by

\begin{equation}\label{ex2}
\left\{
\begin{array}{cc}
{u_1}_t=\nabla(d_1\nabla u_1)+k(x)(bu_3-au_1u_2)&x\in\Omega_1,t>0\\
{u_2}_t=\nabla(d_2\nabla u_2)+k(x)(bu_3-au_1u_2)&x\in\Omega_2,t>0\\
{u_3}_t=\nabla(d_3\nabla u_3)+k(x)(au_1u_2-bu_3)&x\in\Omega_2,t>0\\
\frac{\partial}{\partial\eta}u=0&x\in\partial\Omega_1,t>0\\
\frac{\partial}{\partial\eta}v=\frac{\partial}{\partial\eta}w=0&x\in\partial\Omega_2,t>0\\
u_1=u_{0_1}&x\in\Omega_1,t=0\\
u_2=u_{0_2},u_3=u_{0_3}&x\in\Omega_2,t=0.
\end{array}
\right.
\end{equation}

\begin{figure}[h!]
  \includegraphics[width=\linewidth]{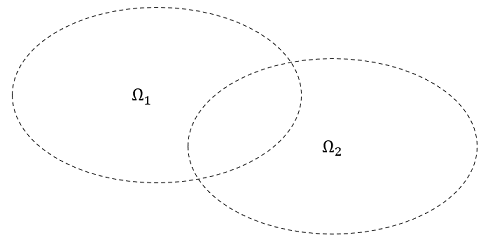}
  \caption{}
  \label{fig:domains}
\end{figure}

\medskip\noindent Here, $d_i$ are positive bounded functions on $\Omega_1$ that are bounded away from $0$, $a,b>0$ and $u_{0_1}$, $u_{0_2}$ and $u_{0_3}$ are nonnegative and bounded. This system has a long history in the setting where $\Omega_1=\Omega_2$, and has appeared in many publications. One of the first was \cite{rothe1984global}, and a multitude of others following. Some of these are cited in \cite{pierre2010global}. 

In the setting when $\Omega_1\ne\Omega_2$, we will see in Section 4 that $u_1$, $u_2$ and $u_3$ are nonnegative. Also, the reaction vector field 
$$
f(x,u)=\left(\begin{array}{c}
k(x)(bu_3-au_1u_2)\\
k(x)(bu_3-au_1u_2)\\
k(x)(au_1u_2-bu_3)
\end{array}
\right)$$
satisfies 
$$f_1(x,u)+f_2(x,u)+2f_3(x,u)=0.$$ 
This guarantees 
$$\|u_1(t,\cdot)\|_{1,\Omega_1}+\|u_2(t,\cdot)\|_{1,\Omega_2}+2\|u_3(t,\cdot)\|_{1,\Omega_2}\le \|u_{0_1}\|_{1,\Omega_1}+ \|u_{0_2}\|_{1,\Omega_2}+2\|u_{0_3}\|_{1,\Omega_2}\text{ for all }t>0.$$
In addition, the component $f_1$ is the only component associated with a species living on all of $\Omega_1$, and it is clearly bounded above by $bu_3$ when $u_i\ge 0$. The two components $f_2$ and $f_3$ corresponding to components associated with species living on all of $\Omega_2$ satisfy 
$$\begin{array}{cc}
f_2(x,u)\le bu_3&\text{for }x\in\Omega_1\cap\Omega_2,u_i\ge 0,\\
f_2(x,u)+f_3(x,u)=0&\text{for }x\in\Omega_1\cap\Omega_2,u_i\ge 0.\\
\end{array}$$
We will see in Section 4 that this structure is sufficient to guaranteed the system (\ref{ex2}) has a unique weak global solution which is sup norm bounded.

\subsection{Notation and Assumptions}
This work focusses on the analysis of reaction-diffusion systems with species on multiple domains. To this end, let $N,n\ge 1$ be integers, and suppose $\Omega_1,...,\Omega_N\subset\R^n$ are bounded domains with smooth boundaries $M_i:= \partial\Omega_i$ for $i=1,...,N$ such that each $\Omega_i$ lies locally on one side of $M_i$. We define $\Omega =\cup _{i=1}^{N}\Omega _{i}$. Each domain $\Omega _{k}$ represents a habitat which houses $n_k$ species of a population having a total of $m$ species. Some of the habitats may overlap, and some may be completely contained in other habitats. We assume there is a mapping $\sigma :\{1,2,...,m\}\rightarrow \{1,2,...,N\}$ which defines species $k$ to be uniquely associated with a habitat $\Omega_{\sigma(k)}$.  Notationally, this results in each species being associated with an appropriate
habitat by partitioning the set $\{1,2,...,m\}$ into $N$ disjoint sets, $O_{1},O_{2},...,O_{N},$ where $i\in O_{j}$ can be interpreted as meaning the $i^{th}$ species is associated with $\Omega _{j}.$ Finally, we denote the population density of species $k$ on $\Omega_{\sigma(k)}$ at time $t\ge 0$ by $u_k(t,\cdot)$.

\medskip We model the interactions of the species $u=\left(u_k\right)_{k=1}^m$ across all habitats via a reaction-diffusion system given by
\begin{equation}
\left\{
\begin{tabular}{ll}
$\frac{\partial}{\partial t}u_k=\nabla\left(d_k(t,x) \nabla u_k\right)+f_k(t,x,u)$ & $t>0,x\in \Omega_{\sigma(k)}\quad k=1,...m$ \\
$\frac{\partial} {\partial \eta_{\sigma(k)}}u_k =0$ & $t>0,x\in M_{\sigma (k)}
\quad k=1,...m$
\\
$u_{k}(0,\cdot )=u_{0_{k}}(\cdot )$ & $t=0,x\in \Omega_{\sigma(k)}
\quad
k=1,...m$
\end{tabular}
\right.   \label{eq1}
\end{equation}
\newline
Here, $u(t,x)=(u_{k}(t,x))_{k=1}^m$ is an unknown vector valued
function. 

\bigskip \noindent\textbf{Assumption (A1):} We assume the structure of the species and habitats described above, and for each $k=1,...,m$, $u_{0_{k}}\in L^\infty(\Omega_{\sigma(k)},\mathbb{R}_+)$, $d_{k}\in L^\infty((0,T), \Omega_{\sigma(k)})$ for each $T>0$, and there exists $\alpha>0$ so that $\alpha\le d_{k}(t,x)$ for all $t>0$ and $x\in\Omega_{\sigma(k)}$. In addition, for each $j=1,...,N$, $\eta_j$ denotes the outward unit normal vector to $\Omega_j$ at a point on $M_j$. For each $k=1,...,m$ we define the $m\times m$ diagonal matrix $Y_k(x)$ for $x\in\Omega_{\sigma(k)}$ so that the $(i,i)$ entry given by the characteristic function $\chi_{\Omega_{\sigma(i)}\cap\Omega_{\sigma(k)}}(x)$, and let $F:\mathbb{R}_+\times \Omega\times \mathbb{R}_+^{m}\rightarrow \mathbb{R}^{m}$ where $F=(F_{k})$, and for each $k=1,...,m$, the function $F_k\in L^{\infty }((0,T)\times\Omega_{\sigma(k)}\times U)$ for bounded subsets $U\subset \mathbb{R}_+^m$ and $T>0$, and $F_k(t,x,u)$ is locally Lipschitz in $u$, uniformly on $(0,T)\times\Omega_{\sigma(k)}$ for each $T>0$. Finally, we define $f=(f_k)$ where $f:\mathbb{R}_+\times \Omega\times \mathbb{R}_+^{m}\to\mathbb{R}^m$ such that 
$$f_k(t,x,u)=\left\{\begin{array}{cc}F_k(t,x,Y_k(x)u),&x\in\Omega_{\sigma(k)}\\0,&\text{otherwise}\end{array}\right.$$ 

\bigskip\noindent We remark that for $k\in\{1,...,\}$, the function $f_k$ has the same qualities as $F_k$, except that for a given $j\in\{1,...,m\}$, $f_k(t,x,u)$ only depends on component $j$ of $u$ if $x\in \Omega_{\sigma(k)}\cap\Omega_{\sigma(j)}$. The extension of $f_k(t,x,u)$ as $0$ outside $\Omega_{\sigma(k)}$ is only done for convenience in development of $L^1$ estimates below.

\bigskip We remark that the homogeneous Neumann boundary conditions listed in (\ref{eq1}) can be replaced with nonhomogeneous boundary conditions. It is also possible to use some ideas from \cite{morgan2021global} to include nondiagonal diffusion, nonlinear diffusion and semilinear boundary conditions. It is also possible to use other simple boundary conditions, including homogeneous Dirichlet boundary conditions. In all cases, it is possible to include convective terms provided $L^1$ apriori estimates can be obtained. The interested reader is referred to \cite{fmty} for additional remarks in the setting of $N=1$, which can be extended with modification to the current setting.

\medskip We are primarily interested in systems which guarantee that solutions to (\ref{eq1}) are componentwise nonnegative, and total population is bounded for finite time. That is, $u_k(t,x)\ge 0$ for each $k=1,...,m$, and there exists $C\in C(\mathbb{R}_+,\mathbb{R}_+)$ such that
\begin{equation}
\sum_{k=1}^m \int_{\Omega_{\sigma(k)}}u_k(t,x)dx\le C(t)
\label{mass}
\end{equation}
for each $t\ge 0$. There has been a wealth of work on systems of the form (\ref{eq1}) when $N=1$. See \cite{desvillettes2007global,souplet2018global,pierre2017dissipative,canizo2014improved,caputo2009global,morgan1990boundedness,hollis1987global,fitzgibbon1997stability, pierre2000blowup, pierre2010global, fellner2021uniform, fellner2020global, desvillettes2007global, morgan2020boundedness, morgan1989global}. The results we present in this work extend some of the work from the setting $N=1$ to $N>1$.

\medskip We start by imposing reasonable conditions on the vector field $f$ to guarantee the nonnegativity of solutions. To this end, we assume 
\begin{equation}
f_k(t,x,u)\ge 0\text{ when } t\ge 0, x\in \Omega_{\sigma(k)}, u\in\mathbb{R}_+^m \text{ and } u_k=0,\text{ for } k=1,...,m.
\label{qp}
\end{equation}
Here, $\mathbb{R}_+^m$ is the set of componentwise nonnegative vectors in $\mathbb{R}^m$. In the setting of $N=1$, condition (\ref{qp}) is typically referred to as a quasi-positivity condition. It is not difficult to prove that solutions to (\ref{eq1}) are componentwise nonnegative regardless of the choice of bounded, componentwise nonnegative initial data if and only if (\ref{qp}) holds. More general information related to nonnegativity of solutions in the case of $N=1$ appears in \cite{CHUEH1}.

\medskip There are many conditions that can result in bounded total population. The one we assume is related to a well known dissipativity condition in the setting $N=1$ that has been used in many of the references listed above (and we especially note \cite{pierre2010global}). The analogous assumption in this setting requires that there exist $b_k>0$ for each $k=1,...,m$, $K_2\ge 0$ and $K_1\in\mathbb{R}$ so that
\begin{equation}
\sum_{k=1}^m b_{k} f_k(t,x,u)\le K_1\sum_{j=1}^m \chi_{\Omega_{\sigma(j)}}(x)u_j+K_2\text{ for }t\ge 0, x\in\Omega \text{ and } u\in\mathbb{R}_+^m,
\label{bal}
\end{equation}
where $\chi_S$ is the characteristic function on the set $S$. It is possible for the constants $K_1$ and $K_2$ in (\ref{bal}) can be replaced by functions depending on $t$ and $x$, and we leave the details to the interested reader. We will see below that this assumption guarantees the estimate given in (\ref{mass}). 

\medskip It is well known in the $N=1$ setting that assumptions (\ref{qp}) and (\ref{bal}) are not sufficient to guarantee the existence of global solutions to (\ref{eq1}) that are sup norm bounded on $(0,T)\times\Omega$ for all $T>0$ (cf \cite{pierre2000blowup,pierre2022blowup}). In fact, when $N=1$, if $\varepsilon>0$, then in the setting when $m=2$ there exist constant diffusion $d_1,d_2>0$, $n\ge 1$, $C>0$, bounded nonnegative initial data, and $f$ satisfying (\ref{qp}), (\ref{bal}) with $|f_k(t,x,u)|\le C(u_1+u_2+1)^{2+\varepsilon}$, such that the solutions to (\ref{eq1}) blow up in the sup norm in finite time \cite{pierre2022blowup}. As a result, we need at least one additional assumption to avoid sup norm blow up in this setting. 

\medskip Recently, in the setting of $N=1$, work in \cite{fmty} proved solutions to (\ref{eq1}) cannot blow up in the sup norm provided there exist $l,C>0$ so that
\begin{equation}
f_k(t,x,u)\le C\left(\sum_{k=1}^m u_k+1\right)^l\text{ for }t>0,x\in\Omega,u\in\mathbb{R}_+^m\nonumber
\end{equation} 
and there exists an $m\times m$ lower triangular matrix $A$ with positive diagonal entries, and a number $1\le r<1+\frac{2}{n}$ so that
\begin{equation}
Af(t,x,u)\le C\vec{1}\left(\sum_{k=1}^m u_k+1\right)^r\text{ for }t>0,x\in\Omega,u\in\mathbb{R}_+^m.\nonumber
\label{int}
\end{equation}
We note that while there is considerable restriction on $r$ in (\ref{int}), there is no restriction on the size of $l$ in (\ref{poly}). In the setting of $N\ge 1$, it is tempting to simply rewrite the assumptions above, but the analysis does not lend itself to the full generalization of the second one.  Instead, we amend the first assumption above to fit our setting, and use more care with the second assumption. To this end, for each $k=1,...,m$, we assume there exist $C,l>0$ (without restriction on size) so that
\begin{equation}
f_k(t,x,u)\le C\left(\sum_{k=1}^m \chi_{\Omega_\sigma(k)}(x)u_k+1\right)^l\text{ for }t>0,x\in\Omega_{\sigma(k)},u\in\mathbb{R}_+^m,
\label{poly}
\end{equation} 
and for each $j=1,...,N$ there is an $n_j\times n_j$ lower triangular matrix $A_j$ with positive entries on the diagonal, and $C,R>0$ with $1\le r<1+\frac{2}{n}$ so that
\begin{equation}
A_jf_{O_j}(t,x,u)\le C\vec{1}\left(\sum_{k=1}^m \chi_{\Omega_\sigma(k)}(x)u_k+1\right)^r\text{ for }t>0,x\in\Omega_j,u\in\mathbb{R}_+^m.
\label{int}
\end{equation}
Here, $f_{O_j}$ denotes the vector whose entries are $f_k$ components of $f$ such that $\Omega_{\sigma(k)}=\Omega_j$. Note that the right hand side of (\ref{int}) includes all components of $u$ whose habitats intersect with $\Omega_j$. 

Note that when components of the vector field are polynomial in nature, the value of $r$ in (\ref{int}) is more restrictive than the inequality indicates. This is because a polynomial bounded above by another polynomial that has a positive integer degree $<M$, tells us the actual bound is of degree $M-1$. So, when $n\ge 2$, the upper bound for $r$ above effectively restricts us to $r=1$, while in the setting of $n=1$, $r$ can be $2$. This does not mean that the reaction terms can only be linear in nature. In deed, we can see in (\ref{eqex1}) that there are quadratic reaction terms, but it is apparent that (\ref{int}) is satisfied with $r=1$.

\medskip The condition in (\ref{int}) has a long history in the setting of $N=1$, and was originally termed an \textit{intermediate sum condition} in \cite{morgan1989global,morgan1990boundedness}. As pointed out in \cite{fmty}, this condition implies a much more general condition that actually leads to the result given in that work, but (\ref{int}) is far easier to recognize in systems, and it occurs naturally as a trade off of higher order terms related to different components. 

\medskip In this work, we extend the results in \cite{fmty} by using (\ref{qp}), (\ref{bal}), (\ref{qp}) and (\ref{int}) to prove that solutions to (\ref{eq1}) cannot blow up in the sup norm in finite time. Section 2 contains some notation, definitions, and the statements of our main results. The proofs are given in Sections 3 and 4, and some examples are stated in Section 5. Finally, we pose an open question in Secton 6.

\section{Statements of Main Results}

Because of the $L^\infty$ nature of the diffusion terms, and the possible abrupt changes in a component $f_k(t,x,u)$ due to dependence on a component $u_j$ for which $\sigma(k)\ne\sigma(j)$ and $\Omega_{\sigma(k)}\cap\Omega_{\sigma(j)}\ne\emptyset$, solutions to (\ref{eq1}) cannot be expected to be classical solutions. As a result, we rely upon the notion of a weak solution.

\begin{definition}\label{def_weak}
	A vector $u = (u_1,\ldots, u_m)$ is called a weak solution to (\ref{eq1}) on $(0,T)$ iff for each $i=1,...,m$,
	\begin{equation*}
		u_i \in C([0,T];L^2(\Omega_{\sigma(i)}))\cap L^2(0,T;H^1(\Omega_{\sigma(i)})), \quad F_i(t,\cdot,u) \in L^2(\Omega_{\sigma(i)}),
	\end{equation*}
	with $u_i(\cdot,0) = u_{0_i}(\cdot)$, and for any test function $\varphi\in L^2(0,T;H^1(\Omega_{\sigma(i)}))$ with \newline$\pa_t\varphi \in L^2(0,T;H^{-1}(\Omega_{\sigma(i)}))$, one has
	\begin{align*}
		&\int_{\Omega_{\sigma(i)}} u_i(t,x)\varphi(t,x)dx\bigg|_{t=0}^{t=T} - \int_0^T\int_{\Omega_{\sigma(i)}} u_i\pa_t\varphi dxdt + \int_0^T\int_{\Omega_{\sigma(i)}} d_i(t,x)\na u_i\cdot \na \varphi dxdt\\
		&=  \int_0^T\int_{\Omega_{\sigma(i)}} f_i(t,x,u)\varphi dxdt.
	\end{align*}
A weak solution to (\ref{eq1}) is a global solution provided it is a weak solution for each $T>0$.
\end{definition}

\medskip\noindent Our main result is stated below.  

\medskip

\begin{theorem}\label{thm1}
	Assume (A1), (\ref{qp}), (\ref{bal}), (\ref{poly}) and (\ref{int}), and that
	\begin{equation}\label{growth}
		1\leq r < 1 + \frac{2}{n}.
	\end{equation}
	Then there exists a unique, componentwise nonnegative, global sup norm bounded weak solution to (\ref{eq1}), i.e. $u_i\in L^\infty_{\rm{loc}}(0,\infty;L^{\infty}(\Omega_{\sigma(i)}))$ for all $i=1,\ldots, m$. Moreover, if $K_1< 0$ or $K_1= K_2 = 0$ in (\ref{bal}), then the solution is bounded uniformly in time. That is,
	\begin{equation}\label{Linftybound}
		{\rm{ess}\sup}_{t\geq 0}\|u_i(t,\cdot)\|_{{\infty},\Omega_{\sigma(i)}} < +\infty, \quad \forall i=1,\ldots, m.
	\end{equation}
\end{theorem}

\medskip

\begin{remark}
	The assumption (\ref{bal}) in Theorem \ref{thm1} is only used to obtain bounds for $u_i(t,\cdot)$ in $L^\infty_{\rm{loc}}(0,\infty;L^{1}(\Omega_{\sigma(i)}))$ for each $i=1,...,m$, and the assumption that $K_1< 0$ or $K_1= K_2 = 0$ in (\ref{bal}) results in bounds for $\|u_i(t,\cdot)\|_{1,\Omega_{\sigma(i)}}$ that are independent of $t$. If these bounds can be obtained by any other means, then Theorem \ref{thm1} remains true without the assumption of (\ref{bal}). In addition, if there exists $a\ge 1$ so that bounds for $u_i(t,\cdot)$ in $L^\infty_{\rm{loc}}(0,\infty;L^{a}(\Omega_{\sigma(i)}))$ can be obtained for each $i=1,...,m$, then the upper bound for $r$ in Theorem \ref{thm1} can be relaxed to
	\begin{equation}
		1\leq r < 1 + \frac{2a}{n}.\nonumber
	\end{equation}
Finally, if bounds for $\|u_i(t,\cdot)\|_{a,\Omega_{\sigma(i)}}$ can be obtained that are independent of $t$, for each $i=1,...,m$, then the uniform sup norm bound in Theorem \ref{thm1} can be obtained. Modifications of the proof given in Section 3 can be employed in accordance with the ideas in \cite{fmty}.
\end{remark}

\medskip We remark that in the case of space dimension $n=1$, we have $2<1+2/n=3$, and as a result, if the components of the reaction vector field $f$ are bounded above by a second degree polynomial, then (\ref{int}) is automatically satisfied. Consequently, we have the simple result below.

\begin{corollary}\label{cor1}
Assume $n=1$, and (A1), (\ref{qp}), (\ref{bal}) and (\ref{poly}) with $l=2$. Then there exists a unique, componentwise nonnegative, global sup norm bounded weak solution to (\ref{eq1}), i.e. $u_i\in L^\infty_{\rm{loc}}(0,\infty;L^{\infty}(\Omega_{\sigma(i)}))$ for all $i=1,\ldots, m$. Moreover, if $K_1< 0$ or $K_1= K_2 = 0$ in (\ref{bal}), then the solution is bounded uniformly in time. That is,
	\begin{equation}\label{Linftybound2}
		{\rm{ess}\sup}_{t\geq 0}\|u_i(t,\cdot)\|_{{\infty},\Omega_{\sigma(i)}} < +\infty, \quad \forall i=1,\ldots, m.
	\end{equation}
\end{corollary}

\medskip We prove Theorem \ref{thm1} in Section 3 by modifying the arguments in \cite{fmty}. That is, for each $\varepsilon>0$, we introduce a system that approximates (\ref{eq1}) and has a unique componentwise nonnegative solution $u^\varepsilon$ that is sup-norm bounded. These approximate systems are constructed in a manner that results in (\ref{qp}), (\ref{bal}), (\ref{poly}) and (\ref{int}) being satisfied in the same manner as (\ref{eq1}). This allows us to utilize the structure guaranteed by (\ref{int}) to employ a modification of the energy functional approach in \cite{fmty} to obtain $L^p(\Omega_{\sigma(k)})$ estimates for $u_k^\varepsilon(t,\cdot)$ for each $k=1,...,m$, $1<p<\infty$ and $t>0$ that are independent of $\varepsilon$. Then, (\ref{poly}) is used, along with results that can be found in either \cite{LSU68} or \cite{nittka2014inhomogeneous} to obtain sup norm estimates for $u^\varepsilon$ that are independent of $\varepsilon>0$. Finally, we pass to the limit as $\varepsilon\to 0$ to obtain convergence to a componentwise nonnegative weak solution to (\ref{eq1}), and uniqueness follows from the local Lipschitz assumption on $f$ in (A1). 

\section{Proof of Theorem \ref{thm1}}

We begin by constructing approximate systems to (\ref{eq1}) similar to the setting where $N=1$ in \cite{fmty}. To this end, for $0<\varepsilon<1$, consider the system
\begin{equation}
\left\{
\begin{tabular}{ll}
$\frac{\partial}{\partial t}u_k^\varepsilon=\nabla\left(d_k(t,x) \nabla u_k^\varepsilon\right)+f_k^\varepsilon(t,x,u_+^\varepsilon)$ & $t>0,x\in \Omega_{\sigma(k)}\quad k=1,...m$ \\
$\frac{\partial} {\partial \eta_{\sigma(k)}}u_k^\varepsilon =0$ & $t>0,x\in M_{\sigma (k)}
\quad k=1,...m$
\\
$u_{k}^\varepsilon(0,\cdot )=u_{0_{k}(\cdot )}$ & $t=0,x\in \Omega_{\sigma(k)}
\quad
k=1,...m$
\end{tabular}
\right.   \label{eq2}
\end{equation}
\medskip \noindent where 
\begin{equation*}
	f_k^\varepsilon(t,x,u_+^\varepsilon):=  f_k(t,x,u_+^\varepsilon)\left(1+\varepsilon\sum_{j=1}^m|f_j(t,x,u_+^\varepsilon)|\right)^{-1},
\end{equation*}
and for $z\in\mathbb{R}^m$, the vector $z_+=((z_k)_+)$ where we define $a_+=\left\{\begin{array}{cc}a&\text{if }a\ge 0\\0&\text{if }a<0\end{array}\right.$ for $a\in\mathbb{R}$. 

\medskip\noindent We remark that the structure of (\ref{eq2}) is similar to (\ref{eq1}), and as a result, we can apply our notion of weak solution to (\ref{eq2}). A slight modification of the arguments in \cite{fmty} allows us to prove that if $T>0$, then there exists a unique weak solution to (\ref{eq2}) on $(0,T)$. Moreover, the construction of the truncated system (\ref{eq2}) allows us to take advantage of the structure of the vector field $f$ that is assumed in Theorem \ref{thm1}. We can easily see that the vector field $f^\varepsilon=(f_k^\varepsilon)$ satisfies (\ref{qp}), (\ref{bal}), (\ref{poly}) and (\ref{int}) in the same manner as $f$, regardless of the choice of $\varepsilon>0$. 

\medskip It is a simple matter to prove $u^\varepsilon$ is componentwise nonnegative. For a given $i=1,...,m$, we choose $\varphi=u_{i,-}^\varepsilon$, where $u_{i,-}^\varepsilon=(-u_i^\varepsilon)_+$. We manipulate the definition of weak solution to show that for all $0<t<T$,
\begin{equation*}
-\frac{1}{2}\int_{\Omega_{\sigma(i)}}|u_{i,-}^\varepsilon(t,x)|^2dx-\alpha\int_0^t\int_{\Omega_{\sigma(i)}}|\nabla u_{i,-}^\varepsilon(s,x)|^2dxds=\int_0^t\int_{\Omega_{\sigma(i)}}u_{i,-}^\varepsilon f_i^\varepsilon(s,x,u_+^\varepsilon(s,x))dxds.
\end{equation*}
\noindent Note that (\ref{qp}) implies the right hand side above is nonnegative. As a result, 
$$\int_{\Omega_{\sigma(i)}}|u_{i,-}^\varepsilon(t,x)|^2dx=0$$
for all $0<t<T$. This implies $u_{i,-}=0$ on $\Omega_{\sigma(i)}$ for all $i=1,...,m$, and consequently, $u^\varepsilon=u_+^\varepsilon$, implying $u^\varepsilon$ is componentwise nonnegative. 

\medskip Now we apply (\ref{bal}) to obtain $L^1$ a priori estimates for $u^\varepsilon$ independent of $0<\varepsilon<1$. We begin by choosing the function $\varphi=1$ in the weak formulation for (\ref{eq2}). This gives
\begin{equation*}
\frac{d}{dt}\int_{\Omega_{\sigma(k)}}u_k^\varepsilon(t,x)dx=\int_{\Omega_{\sigma(k)}}f_k^\varepsilon(t,x,u^\varepsilon(t,x))dx.
\end{equation*}
As a result,  
\begin{equation*}
\frac{d}{dt}\sum_{k=1}^m\int_{\Omega_{\sigma(k)}}b_{k}u_k^\varepsilon(t,x)dx=\sum_{k=1}^m\int_{\Omega_{\sigma(k)}}b_k f_k^\varepsilon(t,x,u)dx=\sum_{k=1}^m\int_{\Omega}b_k f_k^\varepsilon(t,x,u)dx,
\end{equation*}
where the coefficients $b_k$ are associated with (\ref{bal}). As a result, (\ref{bal}) implies 
\begin{equation}\label{nearlygron}
\frac{d}{dt}\sum_{k=1}^m\int_{\Omega_{\sigma(k)}}u_k^\varepsilon(t,x)dx\le K_1\sum_{k=1}^m\int_{\Omega_{\sigma(k)}}u_k^\varepsilon(t,x)dx+K_2.
\end{equation}
Consequently, if $K_1\ne0$, Gronwall's inequality implies
\begin{equation}\label{L11}
\sum_{k=1}^m\|u_k^\varepsilon(t,\cdot)\|_{1,\Omega_{\sigma(k)}}\le \left(\frac{K_2}{K_1}+\sum_{k=1}^m\int_{\Omega_{\sigma(k)}}u_{0_k}(x)dx\right)\exp\left(K_1t\right)-\frac{K_2}{K_1},
\end{equation}
and if $K_1=0$,  Gronwall's inequality implies
\begin{equation}\label{L12}
\sum_{k=1}^m\|u_k^\varepsilon(t,\cdot)\|_{1,\Omega_{\sigma(k)}}\le \sum_{k=1}^m\int_{\Omega_{\sigma(k)}}u_{0_k}(x)dx+K_2t.
\end{equation}
In either case, we have bounds for $\|u_k^\varepsilon(t,\cdot)\|_{1,\Omega_{\sigma(k)}}$ which are independent of $\varepsilon$ for $k=1,...,m$, and the bounds are independent of $t>0$ if either $K_1,K_2=0$ or $K_1<0$.

\medskip Now we use (\ref{int}) to bootstrap the bounds for $\|u_k^\varepsilon(t,\cdot)\|_{1,\Omega_{\sigma(k)}}$ to $\|u_k^\varepsilon(t,\cdot)\|_{p,\Omega_{\sigma(k)}}$ bounds for each $1<p<\infty$. To this end, we build energy functionals in a manner similar to that in \cite{fmty}. Recall that for each $k=1,...,m$, we have $n_k=|O_k|$. Fix $k\in\{1,...,m\}$ and write $\mathbb{Z}_{+}^{n_k}$ as the set of all $n_k$-tuples of non-negative integers. Addition and
	scalar multiplication by non-negative integers of elements in $\mathbb{Z}_{+}^{n_k}$
	is understood in the usual manner. If $\beta=(\beta_{1},...,\beta_{n_k})\in \mathbb{Z}_{+}^{n_k}$ and $p\in \mathbb N\cup \{0\}$,
	then we define $\beta^{p}=((\beta_{1})^{p},...,(\beta_{n_k})^{p})$.
	Also, if $\alpha=(\alpha_{1},...,\alpha_{n_k})\in  \mathbb{Z}_{+}^{n_k}$, then
	we define $|\alpha|=\sum_{i=1}^{n_k}\alpha_{i}$. Finally, if $z=(z_{1},...,z_{n_k})\in \mathbb{R}_{+}^{n_k}$
	and $\alpha=(\alpha_{1},...,\alpha_{n_k})\in \mathbb{Z}_{+}^{n_k}$, then we define
	$z^{\alpha}=z_{1}^{\alpha_{1}}\cdot...\cdot z_{n_k}^{\alpha_{n_k}}$,
	where we interpret $0^{0}$ to be $1$. For simplicity of notation, we momentarily define $v=(u_j^\varepsilon)|_{j\in O_k}$, $g(t,x,u)=(f_j^\varepsilon(t,x,u))|_{j\in O_k}$ and $D(t,x)=(d_j(t,x))|_{j\in O_k}$. Note that each of $v$, $g$ and $D$ have $n_k$ components. For $p\in \mathbb N\cup \{0\}$, we build our $L^p$-energy function of the form
	\begin{equation}\label{Lp}
	\L_{k,p}[v](t) = \int_{\Omega_{\sigma(k)}} \H_p[v](t)dx
	\end{equation}
	where
	\begin{equation}\label{Hp}
	\H_p[v](t) = \sum_{\beta\in \mathbb Z_+^{n_k}, |\beta| = p}\begin{pmatrix}
	p\\ \beta\end{pmatrix}\theta^{\beta^2}v(t)^{\beta},
	\end{equation}
	with
	\begin{equation}\label{our-def}
	\begin{pmatrix}
	p\\ \beta\end{pmatrix}=\frac{p!}{\beta_1!\cdots\beta_{n_k}!},
	\end{equation}
	and $\theta= (\theta_1,\ldots, \theta_{n_k})$ where $\theta_1,...,\theta_{n_k}$ are positive real numbers which will be determined later. For convenience, hereafter we drop the subscript $\beta\in \mathbb Z_+^{n_k}$ in the sum as it should be clear. We note that
	$$\H_0[v](t)=1\text{ and }\H_1[v](t)=\sum_{j\in O_k}\theta_jv_{j}(t).$$
	Also, for a given $p$, $\H_p[v]$ is a general multivariate polynomial of degree $p$ in $v$, and the coefficient defined in (\ref{our-def}) is  the standard multinomial coefficient. Now, suppose $p\ge 2$ is an integer. Proceeding as in \cite{fmty}, let
	 $\L_{k,p}(t):= \L_{k,p}[v](t)$ be defined in \eqref{Lp}.
	Then 
	\begin{align*}
	\frac{d}{dt}\L_{k,p}(t)&=\int_{\Omega_{\sigma(k)}}\sum_{|\beta|=p-1}\left(\begin{array}{c}
	p\\
	\beta
	\end{array}\right)\theta^{\beta^{2}}v(t,x)^{\beta}\sum_{i=1}^{n_k}\theta_{i}^{2\beta_{i}+1}\frac{\partial}{\partial t}v_{i}(t,x)dx\\
	&=\int_{\Omega_{\sigma(k)}}\sum_{|\beta|=p-1}\left(\begin{array}{c}
	p\\
	\beta
	\end{array}\right)\theta^{\beta^{2}}v(t,x)^{\beta}\sum_{i=1}^{n_k}\theta_{i}^{2\beta_{i}+1}\\
	&\qquad\qquad \times \biggl[\na\cdot(D_i(t,x)\na v_{i}(t,x)) +g_{i}(t,x,u^\varepsilon(t,x))\biggr]dx.
	\end{align*}
From \cite{fmty},
	\[
	\int_{\Omega_{\sigma(k)}}\sum_{|\beta|=p-1}\left(\begin{array}{c}
	p\\
	\beta
	\end{array}\right)\theta^{\beta^{2}}v(t,x))^{\beta}\sum_{i=1}^{n_k}\theta_{i}^{2\beta_{i}+1}\na\cdot(D_i(t,x)\na v_i(t,x))dx=I,
	\]
	where 
	\[
	I=-\int_{\Omega_{\sigma(k)}}\sum_{|\beta|=p-2}\left(\begin{array}{c}
	p\\
	\beta
	\end{array}\right)\theta^{\beta^{2}}v(t,x)^{\beta}\sum_{i=1}^{n_k}\sum_{l=1}^{n_k}C_{i,r}(\beta)\left(D_k\nabla v_i(t,x)\right)\cdot\nabla v_l(t,x) dx
	\]
	with
	\[
	C_{i,l}(\beta)=\begin{cases}
	\begin{array}{cc}
	\theta_{i}^{2\beta_{i}+1}\theta_{l}^{2\beta_{l}+1}, & i\ne l,\\
	\theta_{i}^{4\beta_{i}+4}, & i=l.
	\end{array}\end{cases}
	\]
	Then, as in \cite{fmty}, we can can show that for $\theta_i$ sufficiently large, there exists $\alpha_{k,p}>0$ so that
	\begin{align}\label{nearly-there}
	\frac{d}{dt}\L_{k,p}(t)&+\alpha_{k,p}\sum_{i=1}^{n_k} \int_{\Omega_{\sigma(k)}} |\nabla (v_i)^{p/2}(t,x)|^2 dx \nonumber\\
	&\le\int_{\Omega_{\sigma(k)}}\sum_{|\beta|=p-1}\left(\begin{array}{c}
	p\\
	\beta
	\end{array}\right)\theta^{\beta^{2}}v(t,x)^{\beta}\sum_{i=1}^{n_k}\theta_{i}^{2\beta_{i}+1}g_{i}(t,x,u^\varepsilon(t,x))dx.
	\end{align}
	We now look closely at the expression on the right hand side of (\ref{nearly-there}), and in particular the term
	\[
	\sum_{i=1}^{n_k}\theta_{i}^{2\beta_{i}+1}g_{i}(t,x,u^\varepsilon).
	\]
	Note that from (\ref{int}), the definition of the $g_i(t,x,u^\varepsilon)$ and Lemma 2.4 in \cite{fmty}, there exist componentwise increasing functions $h_i:\mathbb{R}^{n_k-i}\to\mathbb{R}_+$ for $i=1,...,n_k-1$ so that if $\gamma_{n_k}>0$ and $\gamma_i\ge h_i(\gamma_{i+1},...,\gamma_{n_k})$ for $i=1,...,n_k-1$ then there exists $K_\gamma>0$ so that
	\[
	\sum_{i=1}^{n_k} \gamma_i g_i(t,x,u^\varepsilon)\le K_\gamma\bra{1+\sum_{i=1}^{m} (\ue_i)^{r}}\text{ for all }(t,x,\ue)\in \mathbb{R}_+\times\Omega_{\sigma(k)}\times\mathbb{R}_+^m.
	\]
	As a result, we can choose $\theta$ so that its components are sufficiently large that the previous positive definiteness condition is satisfied, and 
	\[
	\theta_i\ge h_i(\theta_{i+1}^{2p-1},...,\theta_{n_k}^{2p-1})\text{ for }i=1,...,n_k-1.
	\] 
	Then there exists $K_{\tilde\theta}$ so that for all $\beta\in\mathbb{Z}_+^{n_k}$ with $|\beta|=p-1$, we have
	\[
	\sum_{i=1}^{n_k}\theta_{i}^{2\beta_{i}+1}g_{i}(t,x,\ue(t,x))\le K_{\tilde\theta}\bra{1+\sumi (\ue_i(t,x))^r}\text{ for all }(t,x)\in \mathbb{R}_+\times\Omega_{\sigma(k)}.
	\]
	It follows from this and (\ref{nearly-there}) that there exists $C_p>0$ so that 
	\begin{align}\label{closer-still}
	\frac{d}{dt}\L_{k,p}(t)+\alpha_{k,p}\sum_{i=1}^{n_k} \int_{\Omega_{\sigma(k)}} |\nabla (v_i)^{p/2}(t,x)|^2 dx
	&\le C_p\sum_{j=1}^m\int_{\Omega_{\sigma(j)}}\left(\ue_j(t,x)^{p-1+r}+1\right) dx.
	\end{align}
Now define
$$\L_p(t)=\sum_{k=1}^m\L_{k,p}(t)\text{ and }\alpha_p=\min_{k=1,...,m} \alpha_{k,p}>0.$$
Then from (\ref{closer-still}) and the definition of $v$ for each $k=1,...,m$,
\begin{align}\label{muchcloser-still}
	\frac{d}{dt}\L_{p}(t)+\alpha_{p}\sum_{j=1}^{m} \int_{\Omega_{\sigma(j)}} |\nabla (\ue_j)^{p/2}(t,x)|^2 dx
	&\le mC_p\sum_{j=1}^m\int_{\Omega_{\sigma(j)}}\left(\ue_j(t,x)^{p-1+r}+1\right) dx.
	\end{align}
	Then continuing as in the proof of Theorem 1.1 in \cite{fmty}, there exist $C_p\in C(\mathbb{R}_+,\mathbb{R}_+)$ and $\delta>0$ such that
	\begin{align}\label{muchcloser-still_1}
	\frac{d}{dt}\L_p(t) + \delta \L_p(t)\le C_{p}(t)\quad\forall t>0.
	\end{align}
	Furthermore, $\|C_{p}\|_{\infty,\mathbb{R}_+}<\infty$ if $\|u_k(t,\cdot)\|_{1,\Omega_{\sigma(k)}}$ is bounded independent of $t$ for $k=1,...,m$. Clearly, (\ref{muchcloser-still_1}) allows us to prove there exists $\tilde{C}_p\in C(\mathbb{R}_+,\mathbb{R}_+)$ such that
	\begin{equation*}
		\L_p(t) \le \tilde{C}_{p}(t),\quad \forall t>0,
	\end{equation*}
	with $\|\tilde{C}_{p}\|_{\infty,\mathbb{R}_+}<\infty$ if $\|u_k(t,\cdot)\|_{1,\Omega_{\sigma(k)}}$ is bounded independent of $t$ for $k=1,...,m$. In turn, this allows us to obtain a function $K_p\in C(\mathbb{R}_+,\mathbb{R}_+)$ so that
	\begin{equation*}
		\|\ue_k(t,\cdot)\|_{p,\Omega_{\sigma(k)}} \leq K_{p}(t) \quad \forall t>0, k=1,\ldots, m,
	\end{equation*}
with $\|K_{p}\|_{\infty,\mathbb{R}_+}<\infty$ if $\|u_k(t,\cdot)\|_{1,\Omega_{\sigma(k)}}$ is bounded independent of $t$ for $k=1,...,m$. 

\medskip Finally, we obtain sup norm bounds for $u^\varepsilon$ by using (\ref{poly}) in the same manner as in Proposition 2.1 in \cite{fmty}. Continuing to follow the proof of Theorem 1.1 in \cite{fmty}, we obtain convergence to a solution to (\ref{eq1}), and the remainder of the proof of Theorem \ref{thm1}.

\section{Examples}

In this section, we apply our results to three different example problems given by the system in (\ref{eqex1})-(\ref{eqex1ic}) and (\ref{ex2}), and a model on one dimensional domains that takes illustrates the usefulness of Corollary \ref{cor1}. As we will see, the one dimensional model is a natural follow-up to (\ref{ex2}).

\subsection{Analysis of (\ref{eqex1})-(\ref{eqex1ic})}

To illustrate how Theorem \ref{thm1} applies to (\ref{eqex1})-(\ref{eqex1ic}), we define
$$u=(u_1,u_2,u_3,u_4,u_5,u_6)=(\phi,\psi,\alpha,\beta,v,w),$$ 
$$u=(u_{0_1},u_{0_2},u_{0_3},u_{0_4},u_{0_5},u_{0_6})=(\phi_0,\psi_0,\alpha_0,\beta_0,v_0,w_0)$$
and 
$$
f(t,x,u)=\left(
\begin{array}{c}
-k_1(x)u_1u_4+\lambda_1u_2\chi_{\Omega_1}(x)\\
k_1(x)u_1u_4-\lambda_1u_2\chi_{\Omega_1}(x)\\
-k_2(x)u_3u_2-k_3(x)u_3u_5+\lambda_2u_4\chi_{\Omega_2}(x)\\
k_2(x)u_3u_2+k_3(x)u_3u_5-\lambda_2u_4\chi_{\Omega_2}(x)\\
-k_4(x)u_5u_4\\
k_4(x)u_5u_4-\lambda_3u_6\chi_{\Omega_3}(x)
\end{array}
\right)
$$
for $x\in\Omega=\cup_{i=1}^3\Omega_i$ and $u\in\mathbb{R}_+^6$. Clearly, $f$ satisfies (\ref{qp}) and (\ref{poly}). In addition, (\ref{bal}) is satisfied with $K_1=K_2=0$ since 
$$\sum_{i=1}^6 f_i(x,u)\le 0\quad\text{for }x\in\Omega\text{ and }u\in\mathbb{R}_+^6.$$
Also, $O_1=\{1,2\}$, $O_2=\{3,4\}$ and $O(3)=\{5,6\}$, so
$$f_{O_1}(x,u)=\left(\begin{array}{c}
-k_1(x)u_1u_4+\lambda_1u_2\chi_{\Omega_1}(x)\\
k_1(x)u_1u_4-\lambda_1u_2\chi_{\Omega_1}(x)
\end{array}
\right),
$$
$$
f_{O_2}(x,u)=\left(\begin{array}{c}
-k_2(x)u_3u_2-k_3(x)u_3u_5+\lambda_2u_4\chi_{\Omega_2}(x)\\
k_2(x)u_3u_2+k_3(x)u_3u_5-\lambda_2u_4\chi_{\Omega_2}(x)
\end{array}
\right)
$$
and
$$
f_{O_3}(x,u)=\left(\begin{array}{c}
-k_4(x)u_5u_4\\
k_4(x)u_5u_4-\lambda_3u_6\chi_{\Omega_3}(x)
\end{array}
\right)
$$
for $x\in\Omega$ and $u\in\mathbb{R}_+^6$. So, choosing $A_i=\left(\begin{array}{cc}1&0\\1&1\end{array}\right)$ for $i=1,2,3$ results in (\ref{int}) being satisfied. Therefore, Theorem \ref{thm1} implies (\ref{eqex1})-(\ref{eqex1ic}) has a unique global weak solution, and there exits $C>0$ so that $\|u_i(t,\cdot)\|_{\infty,\Omega_{\sigma(i)}}\le C$ for all $i=1,...,6$ and $t>0$. 

\medskip

\begin{remark}
	If we apply the weak formulation in Definition \ref{def_weak} to $u_5$ and $u_6$ with the test function $\varphi=1$, and sum the results, 
then clearly
$$\int_0^\infty\int_{Omega_3}u_6(t,x)dxdt<\infty.$$
It is possible to use this and further analysis to prove $\|u_6(t,\cdot)\|_{\infty,\Omega_3}\to 0$ as $t\to\infty$. We leave the details and further asymptotic analysis to the interested reader.
\end{remark}

\subsection{Analysis of (\ref{ex2})}

As we pointed out following the statement of (\ref{ex2}), we have
$$
f(x,u)=\left(\begin{array}{c}
k(x)(bu_3-au_1u_2)\\
k(x)(bu_3-au_1u_2)\\
k(x)(au_1u_2-bu_3)
\end{array}
\right)$$
for all $x\in\Omega$ and $u\in\mathbb{R}_+^3$. We can easily see that (\ref{qp}) and (\ref{poly}) are satisfied. In addition, (\ref{bal}) is satisfied  with $K_1=K_2=0$ because 
$$f_1(x,u)+f_2(x,u)+2f_3(x,u)=0$$ 
for all $x\in\Omega$ and $u\in\mathbb{R}_+^3$. In addition, $O_1=\{1\}$ and $O_2=\{2,3\}$, and 
$$f_{O_1}(x,u)=k(x)(bu_3-au_1u_2)\text{ and }f_{O_2}(x,u)=\left(\begin{array}{c}
k(x)(bu_3-au_1u_2)\\
k(x)(au_1u_2-bu_3)
\end{array}\right).$$
Clearly, if $A_1=1$ and
$A_2=\left(\begin{array}{cc}1&0\\1&1\end{array}
\right),$ then (\ref{int}) is satisfied with $r=1$. Therefore, Theorem \ref{thm1} implies (\ref{ex2}) has a unique global weak solution, and there exits $C>0$ so that $\|u_i(t,\cdot)\|_{\infty,\Omega_{\sigma(i)}}\le C$ for all $i=1,2,3$ and $t>0$.

\subsection{A Model in a One Dimensional Setting}

Define $\Omega_1=(0,2)$ and $\Omega_2=(1,3)$, and consider the system given by 

\begin{equation}\label{ex3}
\left\{
\begin{array}{cc}
{u_1}_t=(d_1(t,x)\nabla {u_1}_x)_x+k(x)(u_2^2-u_1u_2)&x\in\Omega_1,t>0\\
{u_2}_t=(d_2(t,x)\nabla {u_2}_x)_x+k(x)(u_1u_2-u_2^2)&x\in\Omega_2,t>0\\
{u_1}_x=0&x\in\{0,2\},t>0\\
{u_2}_x=0&x\in\{1,3\},t>0\\
u_1=u_{0_1}&x\in\Omega_1,t=0\\
u_2=u_{0_2}&x\in\Omega_2,t=0.
\end{array}
\right.
\end{equation}
Here, we assume the functions $d_i$ satisfies (A1) for $i=1,2$, $k(x)$ is the characteristic function on $\Omega_1\cap\Omega_2=(1,2)$, and $u_{0_i}\ge 0$ and bounded for $i=1,2$. If we define
$$f(x,u)=\left(
\begin{array}{c}
k(x)(u_2^2-u_1u_2)\\
k(x)(u_1u_2-u_2^2)
\end{array}
\right),
$$ $O_1={1}$ and $O_2={2}$, then easily (\ref{poly}) and (\ref{qp}) are satisfied, and (\ref{bal}) is satisfied with $K_1=K_2=0$. Furthermore, (\ref{int}) is satisfied with $r=2$, which is admissible in Theorem \ref{thm1} since $n=1$ implies $2<1+r/1=1+2/1=3$. Therefore, Theorem \ref{thm1} implies (\ref{ex3}) has a unique global weak solution, and there exits $C>0$ so that $\|u_i(t,\cdot)\|_{\infty,\Omega_{\sigma(i)}}\le C$ for all $i=1,2$ and $t>0$.

\section{Further Observations and an Open Question}

We point out that that the methods of \cite{fmty} can be  employed as above to extend our results to general advective diffusive operators on each habitat (or domain). Namely, we can obtain unique, globally bounded solutions to when the spatial portion of our differential operators have the form
$$A_i(u_i)= \nabla\cdot\left(D_i(t,x)\nabla u_i+B_i(t,x)u_i\right)$$
for each $i=1,...,m$, where each $D_i\in L^\infty((0,T)\times\Omega_{\sigma(i)},\mathbb{R}^{n\times n})$ is a symmetric positive definite matrix for each $T>0$, and there exists $\delta>0$ so that
$$z^TD_i(t,x)z\ge \delta |z|^2$$
for all $z\in\mathbb{R}^n$. In addition, $B_i\in L^\infty((0,T)\times\Omega_{\sigma(i)},\mathbb{R}^n)$ for each $T>0$. In this case, the homogeneous Neumann boundary conditions in (\ref{eq1}) will be replaced with conditions of the form
$$\left(D_i\nabla u_i)+B_iu_i\right)\cdot \eta=0.$$
It is also possible to obtain results in the setting of quasilinear differential operators. We encourage the interested reader to see \cite{fmty}, and extend the ideas mentioned in that setting. Finally, the boundary conditions listed above can be amended to be nonhomogeneous, and it is also possible to consider homogeneous Dirichlet boundary conditions. The setting of nonhomogeneous Dirichlet boundary conditions presents some problems when it comes to obtaining global existence results from the conditions on the vector field $f$ given earlier. 

\bigskip There many open questions associated with the setting of multiple domains, which arise from the knowledge base associated with the case when $N=1$, even in the setting when the diffusion functions in (\ref{eq1}) are positive constants. We give one of these below. With this in mind, we assume $\Omega_1,\Omega_2\in\mathbb{R}^n$ satisfying the properties listed in Section 1. Also, see figure \ref{fig:domains} in Section 1. For simplicity, assume $d_1,d_2>0$ and consider the system
\begin{equation}\label{eqq}
\left\{
\begin{array}{cc}
u_t=d_1\Delta u+f(x,u,v)&x\in\Omega_1,t>0\\
v_t=d_2\Delta v+g(x,u,v)&x\in\Omega_2,t>0\\
\frac{\partial}{\partial\eta}u=0&x\in M_1,t>0\\
\frac{\partial}{\partial\eta}v=0&x\in M_2,t>0\\
u=u_0&x\in\Omega_1,t=0\\
v=v_0&x\in\Omega_2,t=0
\end{array}
\right.
\end{equation}
\noindent Here, $u_0$ and $v_0$ are bounded nonnegative functions, $f:\Omega_1\times\mathbb{R}_+^2\to\mathbb{R}^2$ and $g:\Omega_2\times\mathbb{R}_+^2\to\mathbb{R}^2$ are locally Lischitz in $u$ and $v$, uniformly in $x$, $f(x,u,v)=0$ for $x\in\Omega_1\backslash\Omega_2$ and $g(x,u,v)=0$ for $x\in\Omega_2\backslash\Omega_1$, $f(x,0,v),g(x,u,0)\ge 0$ for $u,v\ge 0$ and $x\in\Omega_1\cap\Omega_2$, and $f(x,u,v)+g(x,u,v)=0$ for $x\in\Omega_1\cap\Omega_2$ and $u,v\ge 0$. 

\medskip\noindent\textbf{Open Question:} In the setting where $\Omega_1=\Omega_2$, if $f$ and $g$ are smooth, and both are bounded in absolute value by a quadratic polynomial in $u$ and $v$, the results in \cite{fellner2020global,fellner2021uniform} guarantee global existence and uniform sup norm bounds for solutions to (\ref{eqq}). This is an open problem in the setting when $\Omega_1\ne\Omega_1\cap\Omega_2\ne\emptyset$.

\newcommand{\etalchar}[1]{$^{#1}$}

\end{document}